\documentclass[12pt]{article}

\usepackage{amsmath,amssymb}
\usepackage{graphicx,subfig}
\usepackage{parskip}
\usepackage{hyperref}
\usepackage{vmargin}
\setpapersize{USletter}
\setmarginsrb{2.5cm}{1.5cm}{2.5cm}{1.5cm}{0cm}{0cm}{1cm}{1cm}

\begin{document}

\title{Correction to ``Parametric Resonance in Immersed Elastic
  Boundaries''\\
  {\footnotesize (SIAM. J. Appl. Math. \underline{65}(2):494-520, 2004)}
}
\author{William~Ko and John~M. Stockie\\
  Department of Mathematics\\
  Simon Fraser University}

\date{July 19, 2012}

\maketitle

\renewcommand{\thefigure}{\Alph{figure}}

\begin{abstract}
  This note is a correction to a paper of Cortez, Peskin, Stockie \&\
  Varela [\emph{SIAM J. Appl. Math.}, 65(2):494-520, 2004], who studied
  the stability of a parametrically-forced, circular, elastic fiber
  immersed in an incompressible fluid in 2D, and showed the existence of
  parametric resonance.  The results were represented as plots that
  separate parameter space into regions where the solution is either
  stable or unstable.  We uncovered two errors in the paper: the first
  was in the derivation of the eigenvalue problem, and the second was in
  the code to used to calculate the stability contours.
\end{abstract}

An error in the derivation was found in Appendix~B of \cite{Cortez04},
where the gradient of the Dirac delta function was improperly
transformed from Cartesian to polar coordinates.  This led to an error
in the perturbation expansion of the immersed boundary forcing term.
The correct transformation is~\cite{hassani-ch5-2009}
\begin{gather*}
  \nabla\delta(\boldsymbol{x}-\boldsymbol{X}^{(0)}) 
  = \frac{1}{r}\nabla[ \delta(r-1) \delta(\theta-s) ],
\end{gather*}
which changes the second equation in {\scshape Claim}~1 to
\begin{gather*}
  (\boldsymbol{\hat{z}}\cdot\nabla\times\boldsymbol{f}^{(1)})
  = K(t)(X_{ss}^\theta + X_s^r)\left(\frac{\delta_r(r-1)}{r}\right)
  - K(t)(X_{sss}^r - X_{ss}^\theta) \frac{\delta(r-1)}{r},
\end{gather*}
where subscripts denote partial derivatives.  This in turn modifies
equations (6.1)--(6.4) in \cite{Cortez04} for the eigenvalue problem.
For the case when $(\gamma + i n) \ne 0$, the equations are
\begin{align*}
  0 =\; & i\left\{ \phi\Omega_n^3 
    \left[ \frac{H_p(i\Omega_n)}{H_{p-1}(i\Omega_n)} 
      - \frac{J_p(i\Omega_n)}{J_{p+1}(i\Omega_n)}
    \right] + ip
  \right\} X_n^r \\
  \nonumber
  & + \left\{ \phi\Omega_n^3
    \left[ \frac{H_p(i\Omega_n)}{H_{p-1}(i\Omega_n)}
      + \frac{J_p(i\Omega_n)}{J_{p+1}(i\Omega_n)}
    \right] - ip^2
  \right\} X_n^\theta \\
  \nonumber
  & + i\tau p (X_{n-1}^r - X_{n+1}^r) - \tau p^2 (X_{n-1}^\theta - X_{n+1}^\theta), \\
  0 = \; & i\left\{ \phi\Omega_n^4 
    \left[2 - \frac{H_{p+1}(i\Omega_n)}{H_{p-1}(i\Omega_n)} 
      - \frac{J_{p-1}(i\Omega_n)}{J_{p+1}(i\Omega_n)}
    \right] + 2p^3
  \right\} X_n^r \\
  \nonumber
  & - \left\{ \phi\Omega_n^4
    \left[ \frac{H_{p+1}(i\Omega_n)}{H_{p-1}(i\Omega_n)}
      - \frac{J_{p-1}(i\Omega_n)}{J_{p+1}(i\Omega_n)}
    \right] + 2p^2
  \right\} X_n^\theta \\
  \nonumber
  & + 2\tau p^3 (X_{n-1}^r - X_{n+1}^r) + 2\tau ip^2 (X_{n-1}^\theta - X_{n+1}^\theta),
\end{align*}
while when $(\gamma + i n) = 0$ we have
\begin{align*}
  0 &= ip^2 (X_n^\theta - i\tau X_{n-1}^\theta + i\tau X_{n+1}^\theta)
  + p (X_n^r - i\tau X_{n-1}^r + i\tau X_{n+1}^r), \\
  0 &= p^2 (X_n^\theta - i\tau X_{n-1}^\theta + i\tau X_{n+1}^\theta)
  + ip^3 (X_n^r - i\tau X_{n-1}^r + i\tau X_{n+1}^r).
\end{align*}

Note that in each pair of equations, the correction only affects the
second equation while the first remains unchanged.  These results were
derived with the Navier-Stokes equations written in stream
function--vorticity variables.  The derivation was repeated in terms of
primitive variables, and identical results were obtained.

\begin{figure}[ht!]
  \centering
  \fbox{\includegraphics[width=0.9\textwidth]{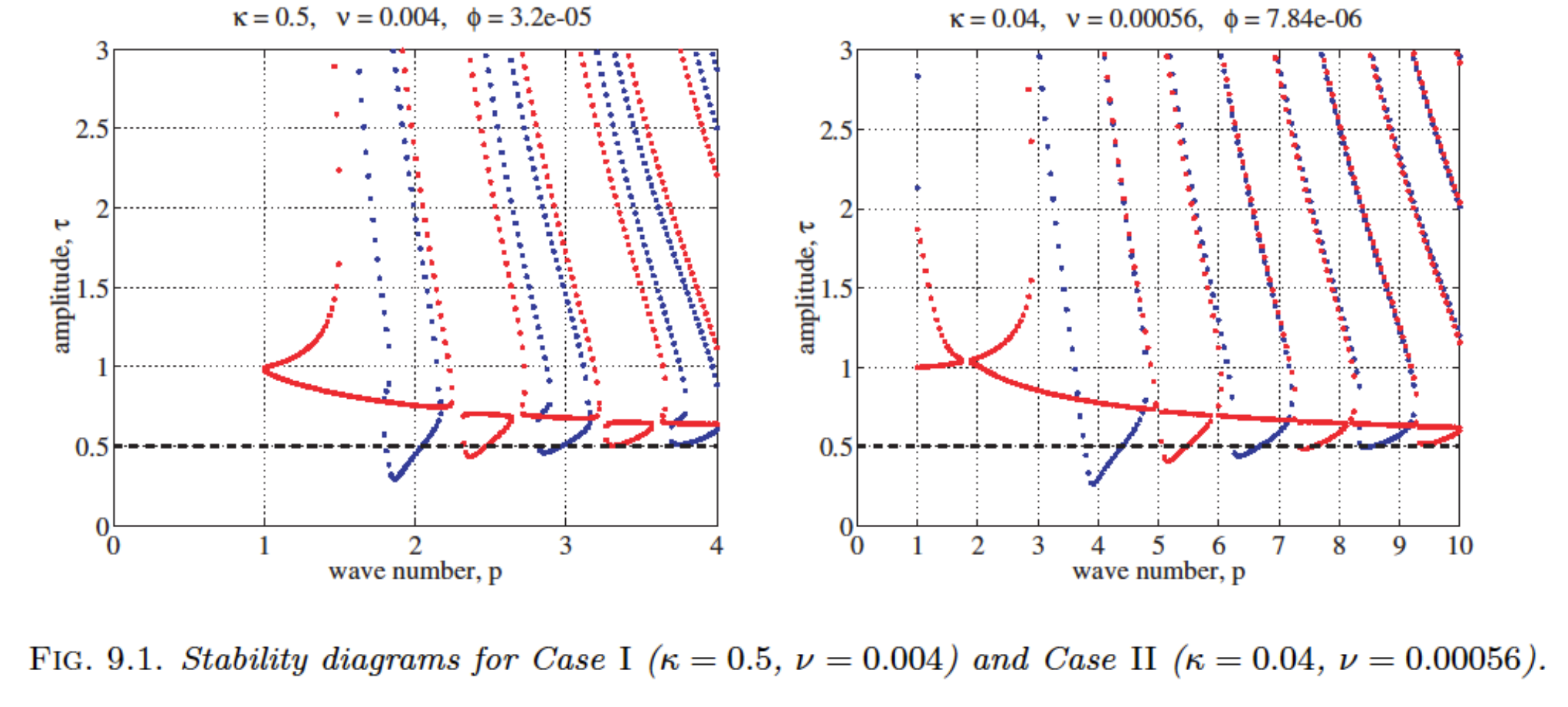}}
  \caption{Original stability contours for cases I and II, reproduced
    from \cite{Cortez04}.}
  \label{fig:fig9p1}
\end{figure}

\begin{figure}[ht!]
  \centering
  \includegraphics[width=0.9\textwidth]{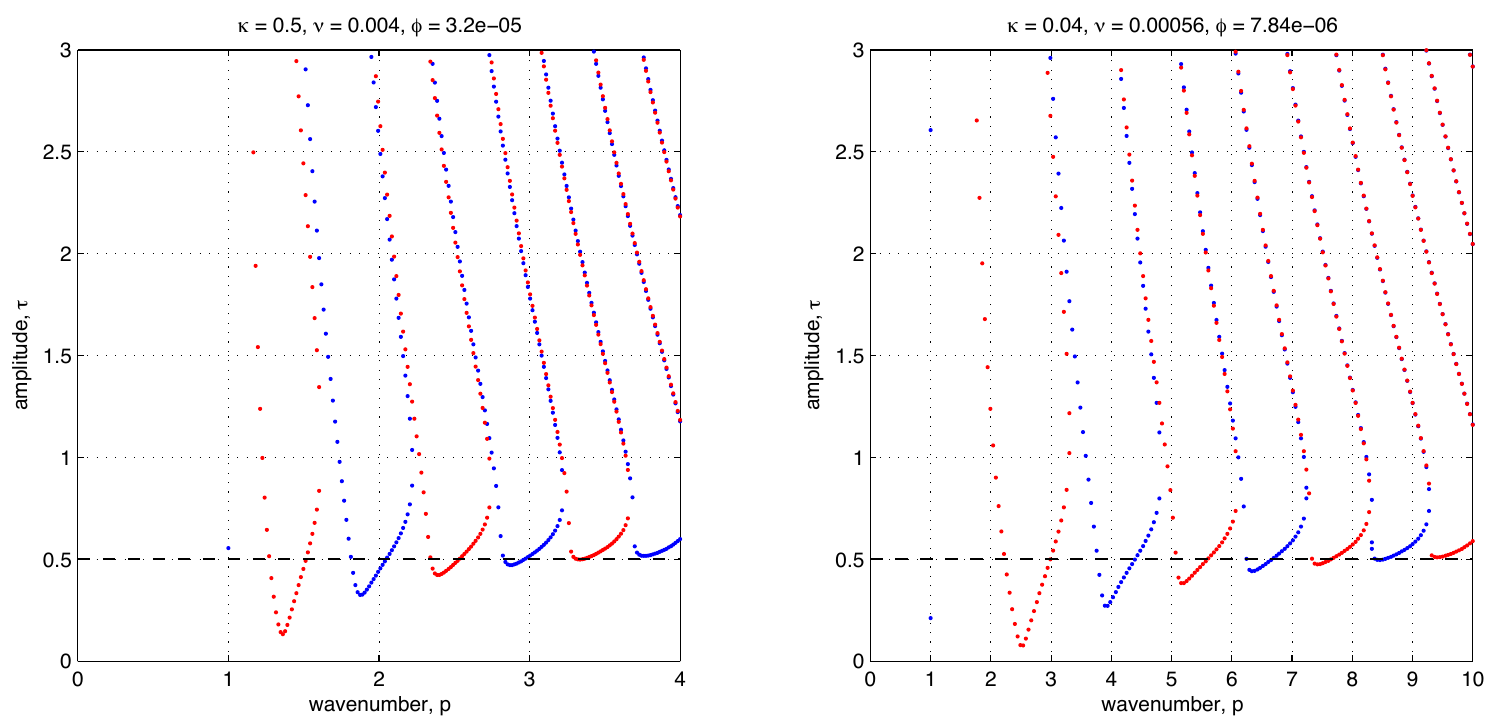}
  \caption{Corrected stability contours for cases I and II (compare to
    Figure~\ref{fig:fig9p1}).}
  \label{fig:case1-2}
\end{figure}

The second error was in the Matlab code used to generate the stability
plots, which affected the stability plots in three ways:
\begin{itemize}
\item One solution mode was missing for the lowest wavenumber
  ``subharmonic'' case.
\item Spurious points were obtained in the region where the subharmonic
  mode was supposed to have been located (i.e., for small values of
  wavenumber $p$), as well as in the upper right reaches of the plot
  (for large values of both wavenumber and amplitude, see especially the
  rightmost plot in Figure~8.2).  There was also a spurious line of points
  running just above the line $\tau=0.5$.
\item The borders of the stability regions were in some cases
  noticeably deformed.
\end{itemize}
All of these discrepancies can be seen by comparing the original
Figure~\ref{fig:fig9p1} (which repeats the stability plots for cases I
and II from Figure~9.1 in \cite{Cortez04}) with the new
Figure~\ref{fig:case1-2} (which shows the corrected plots generated
using the new equations).  Both of these plots correspond to the
following parameter values:
\begin{center}
  \begin{tabular}{ll}
    Case~I & Case~II  \\ \hline
    $\kappa = 0.5$            & $\kappa$ = 0.04 \\
    $\nu = 0.004$             & $\nu = 0.00056$ \\
    $\phi = 3.2\times10^{-5}$ & $\phi = 7.84\times10^{-6}$
  \end{tabular}
\end{center}

Each point depicted in these plots represents one solution to the
eigenvalue problem formulated in Eq.~(8.1) of \cite{Cortez04}, and the
stability boundaries are traced out by varying the angular wavenumber
($p$) of the circular fiber.  The blue points denote harmonic modes and
the red points denote subharmonic modes.  The only ``physical''
instabilities are those corresponding to integer values of wavenumber
$p$ that lie inside the stability ``tongues'' and satisfy $0\leqslant
\tau \leqslant \frac{1}{2}$.  With this in mind, the lowest wavenumber
unstable mode in case~I corresponds to $p=2$, while that in case~II has
$p=4$.

\begin{figure}[tbp]
  \centering
  \fbox{\includegraphics[width=0.9\textwidth]{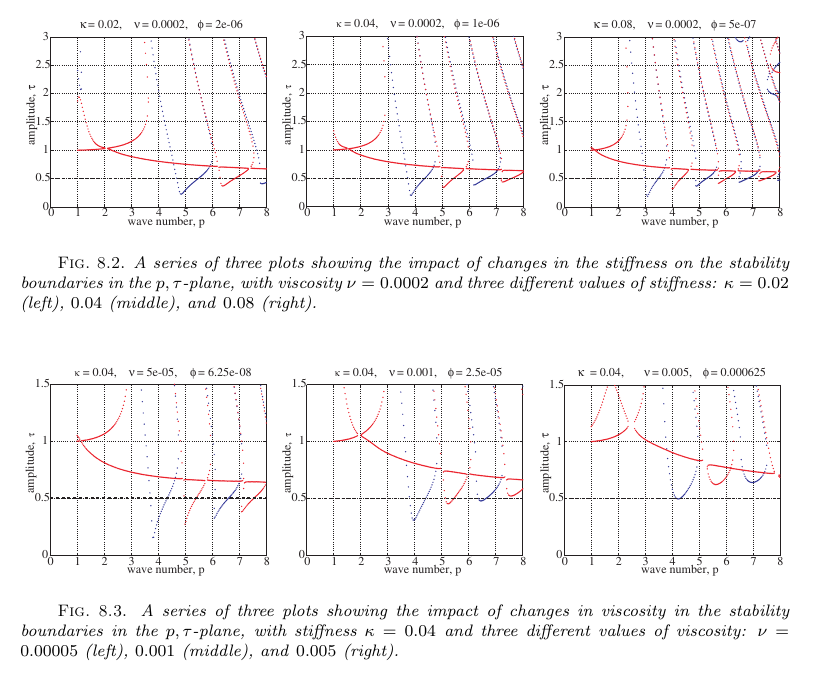}}
  \caption{Original stability plots, reproduced 
    from \cite{Cortez04}.} 
  \label{fig:cortez}
\end{figure}

\begin{figure}[tbp]
  \centering
  \includegraphics[width=0.9\textwidth]{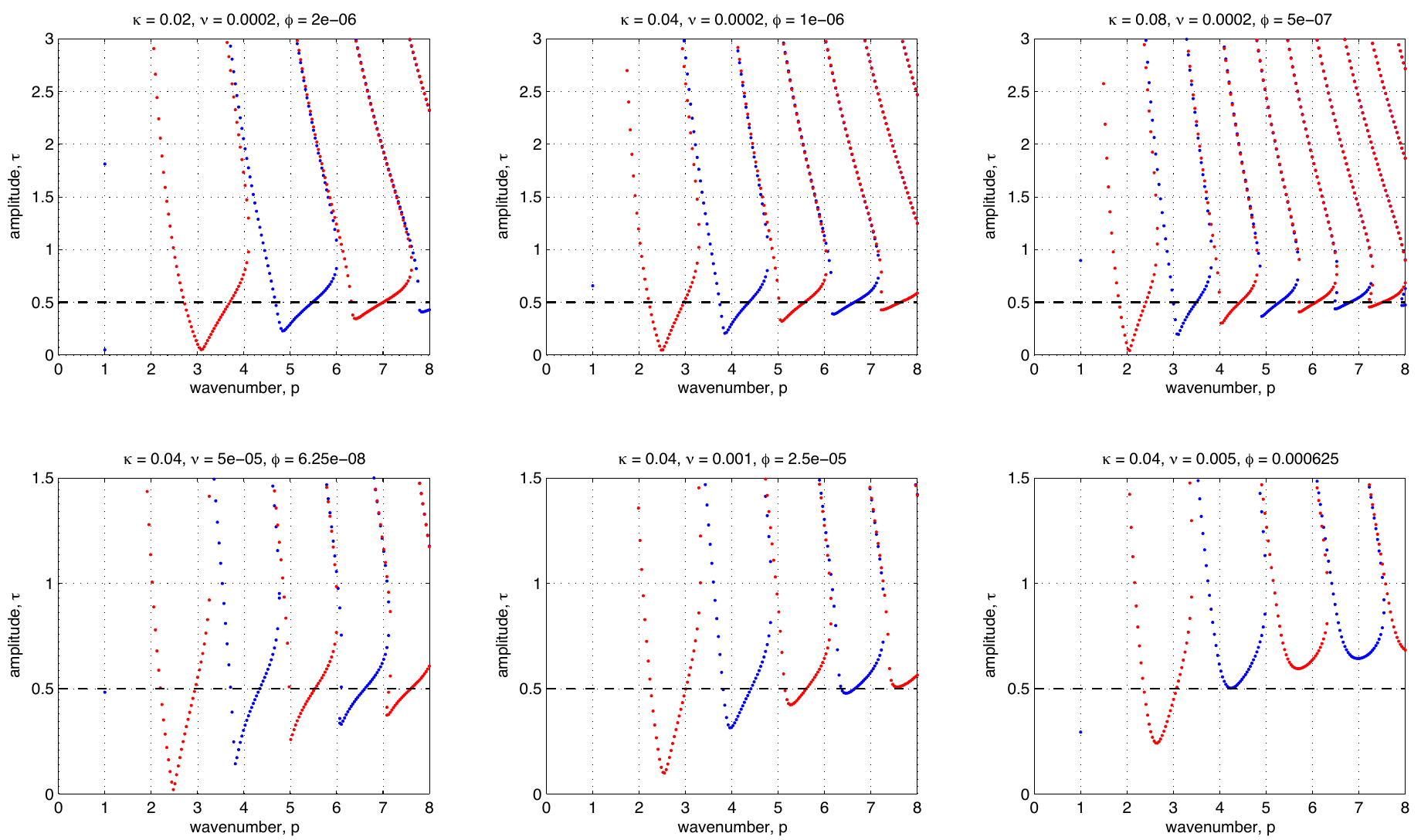}
  \caption{Corrected stability plots (compare to
    Figure~\ref{fig:cortez}).}
  \label{fig:contours}
\end{figure}

Other than the discrepancies noted above, the original stability plots
from \cite{Cortez04} are still qualitatively correct and the parametric
resonances identified in the original paper are also true instabilities
(as verified in numerical simulations).  Nevertheless, we have performed
new calculations with the corrected Matlab code that demonstrate the
existence of resonant subharmonic modes that were not identified in the
original paper.  New resonant modes occur only when the missing
subharmonic region overlaps with an integer value of $p$, which is not
the case in the two examples plotted in Figures~\ref{fig:fig9p1} and
\ref{fig:case1-2}.  The original stability plots from Figures~8.2
and~8.3 are reproduced in Figure~\ref{fig:cortez} while our corrected
plots are given in Figure~\ref{fig:contours}.  In all cases, the
spurious points are eliminated and a new subharmonic mode appears.

\begin{figure}[tbp]
  \centering
  \subfloat{\includegraphics[width=0.45\textwidth]{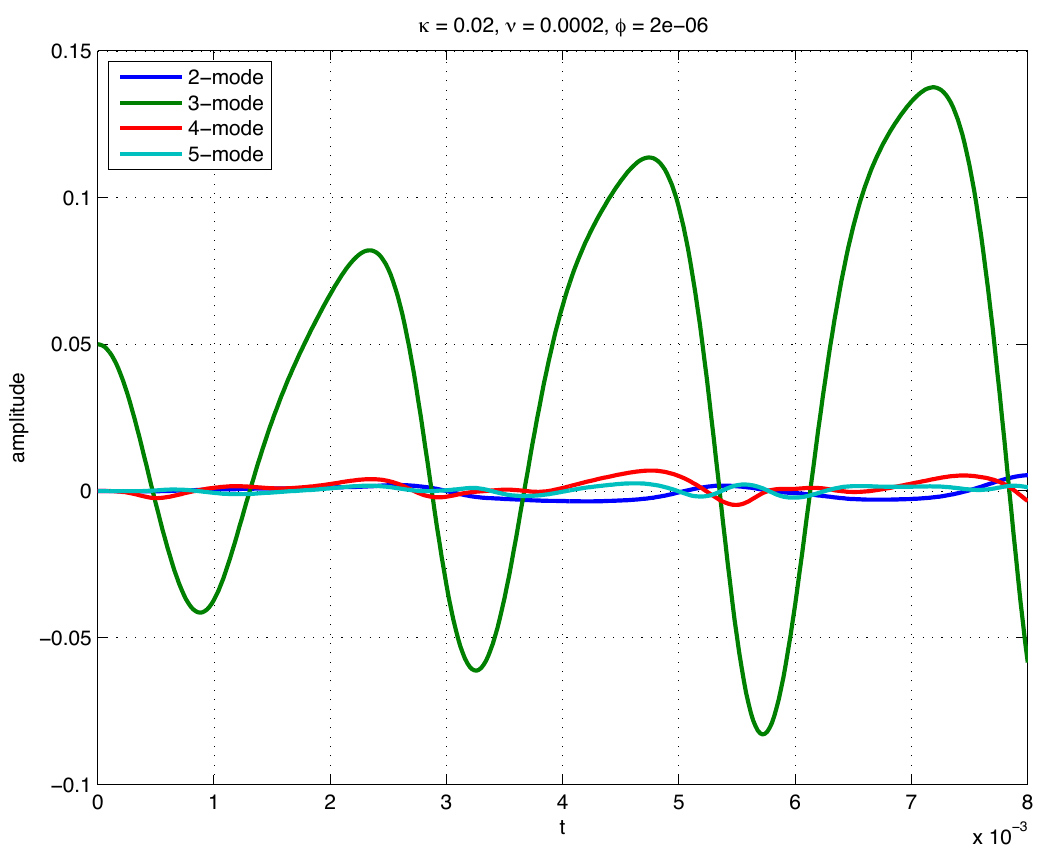}}
  \subfloat{\includegraphics[width=0.45\textwidth]{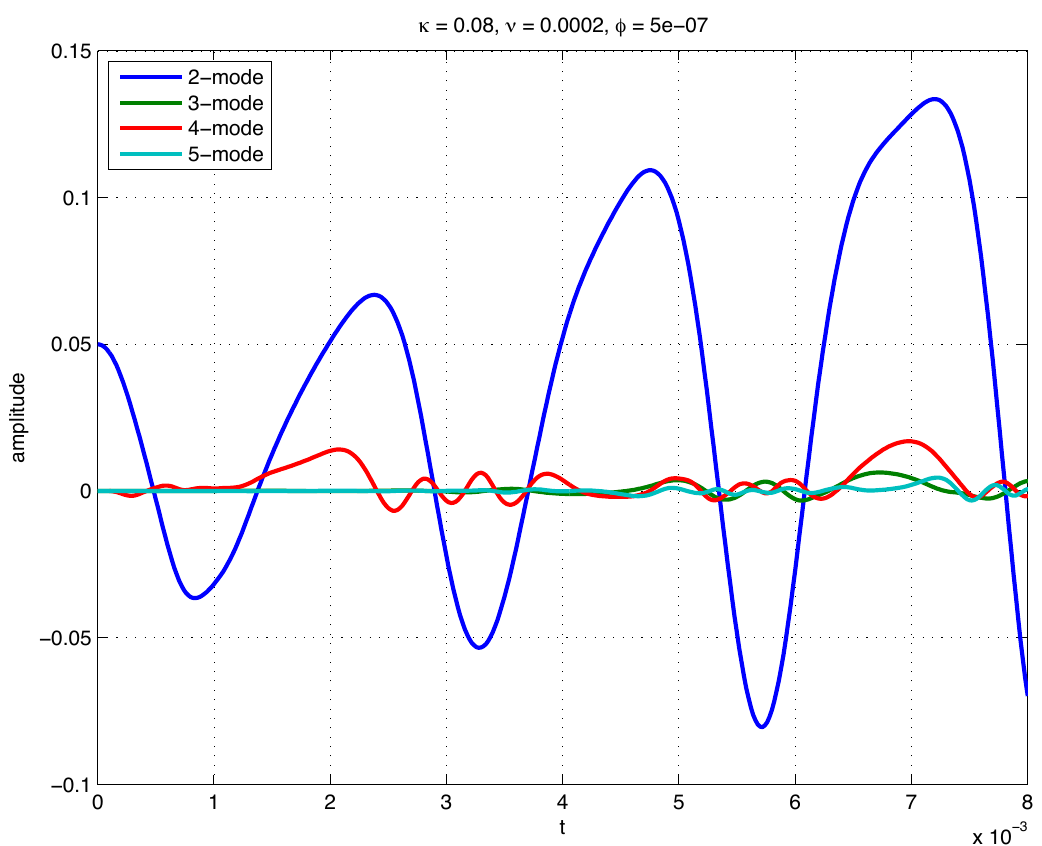}}
  \caption{Amplitudes of the various $p$-modes from numerical
    simulations of Cases~III (left, $p=3$ unstable) and IV (right, $p=2$
    unstable).}
  \label{fig:case3-4}
\end{figure}

To illustrate the presence of the missing unstable modes, we choose
the following two sets of parameters listed below
\begin{center}
  \begin{tabular}{ll}
    Case~III & Case~IV  \\ \hline
    $\kappa = 0.02$         & $\kappa = 0.08$ \\
    $\nu = 0.0002$          & $\nu = 0.0002$ \\
    $\phi = 2\times10^{-6}$ & $\phi = 5\times10^{-7}$
  \end{tabular}
\end{center}
which correspond (respectively) to the left- and right-most plots in
Figures~8.2-8.3.  Numerical simulations of the full immersed boundary
equations have been performed for Case III (with $p=3$) and Case IV
(with $p=2$) to verify that these parametric resonances are actual
instabilities.  Following the approach used in \cite{Cortez04}, we
initialize the fiber using the perturbed circular shape $r =
R(1+0.05\cos(p\theta))$, where $p$ is chosen equal to the resonant
wavenumber.  Figure~\ref{fig:case3-4} shows the amplitude of various
$p$-modes in the two simulations, from which it is clear that the
resonant mode is excited and its amplitude grows in time, as would be
expected for the resonant case.  The amplitude of the other modes
remains small.

\begin{figure}[tbp]
  \centering
  \subfloat{\includegraphics[width=0.45\textwidth]{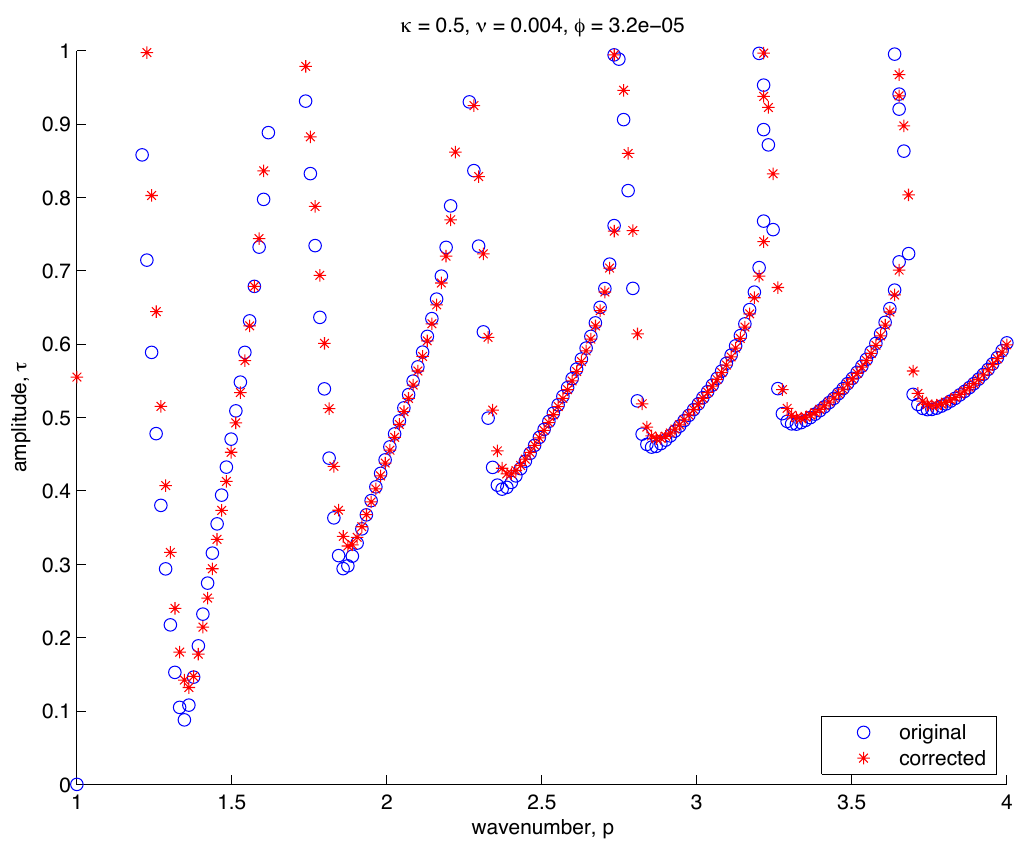}}
  \subfloat{\includegraphics[width=0.45\textwidth]{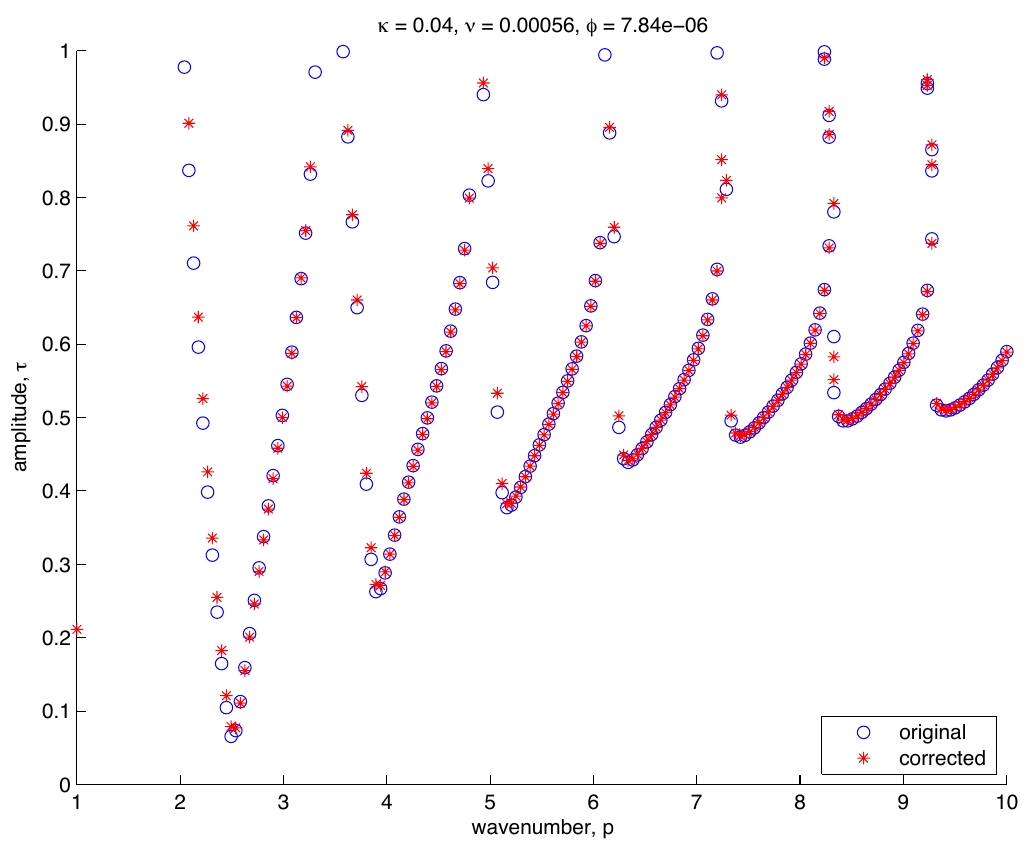}}
  \caption{Comparison of the eigenvalue plots for the original (blue)
    and corrected (red) equations in Cases~I (left) and II (right).}
  \label{fig:compare}
\end{figure}

Other than the presence of these missing lowest-wavenumber subharmonic
modes, the effect on the stability contours is relatively small for the
parameter values we considered.  Figure~\ref{fig:compare} draws a direct
comparison between the stability regions for cases~I and II obtained
with the original equations and with the corrected equations.  Clearly,
the correction results in only a small deformation in the contours and
thus has only minimal impact on the solution.  All other cases
considered in~\cite{Cortez04} yield similar results.

\bibliographystyle{abbrv}
\bibliography{biblio}

\begin{thebibliography}{1}

\bibitem{Cortez04}
R.~Cortez, C.~S. Peskin, J.~M. Stockie, and D.~Varela.
\newblock Parametric resonance in immersed elastic boundaries.
\newblock {\em SIAM Journal on Applied Mathematics}, 65(2):494--520, 2004.

\bibitem{hassani-ch5-2009}
S.~Hassani.
\newblock Dirac delta function.
\newblock In {\em Mathematical Methods for Students of Physics and Related
  Fields}, chapter~5, pages 139--170. Springer, New York, 2009.
\newblock \url{http://dx.doi.org/10.1007/978-0-387-09504-2_5}.

\end{thebibliography}

\end{document}